\documentclass[a4paper, 11pt]{amsart}
\usepackage[utf8]{inputenc}
\usepackage[T1]{fontenc}
\usepackage[english]{babel}
\usepackage{graphicx}
\usepackage{stmaryrd}
\usepackage{tikz}
\usepackage{tikz-cd} 
\usepackage[all]{xy}
\usepackage{comment}
\usepackage{bm}
\usepackage{amsmath,amsfonts,amssymb}
\usepackage[a4paper]{geometry}
\usepackage{amsthm}
\usepackage{float}
\usepackage{enumitem}
\usepackage{hyperref}
\usepackage{xcolor}
\hypersetup{
    colorlinks=true,
    linkcolor={blue},
    citecolor={blue},
    urlcolor={blue}}
\newtheorem{te}{Theorem}[section]

\newtheorem{prop}[te]{Proposition}

\newtheorem{co}[te]{Corollary}

\newtheorem{lemme}[te]{Lemma}
\theoremstyle{definition}

\newtheorem{de}[te]{Definition}

\theoremstyle{remark}
\newtheorem{rque}[te]{Remark}
\newenvironment{proofL}{\noindent\textit{Proof of Lemma \ref{5inv}.~}}{\hfill$\square$\bigbreak} 
\newenvironment{proofTH1}{\noindent\textit{Proof of Theorem \ref{1}.~}}{\hfill$\square$\bigbreak} 
\newenvironment{proofPROP}{\noindent\textit{Proof of Proposition \ref{256}.~}}{\hfill$\square$\bigbreak}

\newlength{\plarg}
\usetikzlibrary{cd}
\calclayout
\title{Elementary subgroups of virtually free groups}
\author{Simon André}
\begin{document}

\begin{minipage}{\linewidth}
	\begin{abstract}
We give a description of elementary subgroups (in the sense of first-order logic) of finitely generated virtually free groups. In particular, we recover the fact that elementary subgroups of finitely generated free groups are free factors. Moreover, one gives an algorithm that takes as input a finite presentation of a virtually free group $G$ and a finite subset $X$ of $G$, and decides if the subgroup of $G$ generated by $X$ is $\exists\forall\exists$-elementary. One also proves that every elementary embedding of an equationally noetherian group into itself is an automorphism.
	\end{abstract}
	\maketitle
\end{minipage}


\thispagestyle{empty}

\section{Introduction}



A morphism $\varphi : H \rightarrow G$ between two groups $H$ and $G$ is said to be \emph{elementary} if the following condition holds: for every first-order formula $\theta(x_1,\ldots,x_k)$ with $k$ free variables in the language of groups, and for every $k$-tuple $(h_1,\ldots,h_k)\in H^k$, the statement $\theta(h_1,\ldots,h_k)$ is true in $H$ if and only if the statement $\theta(\varphi(h_1),\ldots,\varphi(h_k))$ is true in $G$. In particular, $\varphi$ is injective. When $H$ is a subgroup of $G$ and $\varphi$ is the inclusion of $H$ into $G$, one says that $H$ is an \emph{elementary subgroup} of $G$. If one only considers a certain fragment $\mathcal{F}$ of the set of first-order formulas (for instance the set of $\forall\exists$-formulas, $\exists\forall\exists$-formulas or $\exists^{+}$-formulas, see paragraph \ref{logic} for definitions), one says that $\varphi$ (or $H$) is \emph{$\mathcal{F}$-elementary}.

It was proved by Sela in \cite{Sel06} and by Kharlampovich and Myasnikov in \cite{KM06} that any free factor of a non-abelian finitely generated free group is elementary. Later, Perin proved that the converse holds: if $H$ is an elementary subgroup of $F_n$, then $F_n$ splits as a free product $F_n=H\ast H'$ (see \cite{Per11}). Recently, Perin gave another proof of this result (see \cite{Per19}). More generally, Sela \cite{Sel09} and Perin \cite{Per11} described elementary subgroups of torsion-free hyperbolic groups.

Our main theorem provides a characterization of $\exists\forall\exists$-elementary subgroups of virtually free groups. Recall that a group is said to be virtually free if it has a free subgroup of finite index. In what follows, all virtually free groups are assumed to be finitely generated and non virtually cyclic (here, and in the remainder of this paper, virtually cyclic means finite or virtually $\mathbb{Z}$). In \cite{And19}, we classified virtually free groups up to $\forall\exists$-elementary equivalence, i.e.\ we gave necessary and sufficient conditions for two virtually free groups $G$ and $G'$ to have the same $\forall\exists$-theory. In this context, we introduced Definition \ref{legal} below. Recall that a non virtually cyclic subgroup $G'$ of a hyperbolic group $G$ normalizes a unique maximal finite subgroup of $G$, denoted by $E_G(G')$ (see \cite{Ol93} Proposition 1).

\begin{de}[\emph{Legal large extension}]\label{legal}Let $G$ be a non virtually cyclic hyperbolic group, and let $H$ be a subgroup of $G$. One says that $G$ is a \emph{legal large extension} of $H$ if there exists a finite subgroup $C$ of $H$ such that the normalizer $N_H(C)$ of $C$ is non virtually cyclic, the finite group $E_H(N_H(C))$ is equal to $C$, and $G$ admits the following presentation: \[G=\langle H,t \ \vert \ [t,c]=1, \ \forall c\in C\rangle.\]More generally, one says that $G$ is a \emph{multiple legal large extension} of $H$ if there exists a finite sequence of subgroups $H=G_0\subset G_1\subset \cdots \subset G_n=G$ such that $G_{i+1}$ is a legal large extension of $G_i$ for every integer $0\leq i\leq n-1$, with $n\geq 1$.
\end{de}

In terms of graphs of groups, $G$ is a multiple legal large extension of $H$ if it splits as a finite graph of groups over finite groups, whose underlying graph is a rose and whose central vertex group is $H$, with additional assumptions on the edge groups. The prototypical example of a multiple legal large extension is given by the splitting of the free group $G=F_k$ of rank $k\geq 3$ as $F_k=\langle F_2,t_1,\ldots,t_{k-2} \ \vert \ \varnothing\rangle$. In this example, $H$ is the free group $F_2$.

In \cite{And19}, we proved the following result (see Theorem 1.10).

\begin{te}\label{legalteplus}Let $G$ be a non virtually cyclic hyperbolic group, and let $H$ be a subgroup of $G$. If $G$ is a multiple legal large extension of $H$, then $H$ is $\exists\forall\exists$-elementary.\end{te}

\begin{rque}We conjectured in \cite{And19} that $H$ is elementary.
\end{rque}

We shall prove that the converse of Theorem \ref{legalteplus} holds, provided that $G$ is a virtually free group.

\begin{te}\label{1}Let $G$ be a virtually free group, and let $H$ be a proper subgroup of $G$. If $H$ is $\forall\exists$-elementary, then $G$ is a multiple legal large extension of $H$.\end{te}

\begin{rque}In particular, Theorem \ref{1} recovers the result proved by Perin in \cite{Per11}: an elementary subgroup of a free group is a free factor.\end{rque}

\begin{rque}In our classification of virtually free groups up to $\forall\exists$-elementary equivalence (see \cite{And19}), another kind of extension, called legal small extension, plays an important role. Theorem \ref{1} above shows that if $G$ is a non-trivial legal small extension of $H$, then $H$ is not an elementary subgroup of $G$. See also \cite{And19} Remark 1.15.\end{rque}

Putting together Theorem \ref{1} and Theorem \ref{legalteplus}, we get the following result.

\begin{te}\label{2dejapris}Let $G$ be a virtually free group, and let $H$ be a proper subgroup of $G$. The following assertions are equivalent:
\begin{enumerate}
\item $H$ is $\exists\forall\exists$-elementary;
\item $H$ is $\forall\exists$-elementary;
\item $G$ is a multiple legal large extension of $H$.
\end{enumerate}
\end{te}

In addition, one gives an algorithm that decides whether or not a finitely generated subgroup of a virtually free group is $\exists\forall\exists$-elementary.

\begin{te}There is an algorithm that, given a finite presentation of a virtually free group $G$ and a finite subset $X \subset G$, outputs `Yes' if the subgroup of $G$ generated by $X$ is $\exists\forall\exists$-elementary, and `No' otherwise.\end{te}

\begin{rque}Note that any $\forall\exists$-elementary subgroup of a virtually free group is finitely generated, as a consequence of Theorem \ref{1}.
\end{rque}

Recall that every virtually free group $G$ splits as a finite graph of finite groups (which is not unique), called a Stallings splitting of $G$. The following result is an immediate consequence of Theorem \ref{1}.

\begin{co}Let $G$ be a virtually free group. If the underlying graph of a Stallings splitting of $G$ is a tree, then $G$ has no proper elementary subgroup.\end{co}

For instance, the virtually free group $\mathrm{SL}_2(\mathbb{Z})$, which is isomorphic to $\mathbb{Z}/4\mathbb{Z}\ast_{\mathbb{Z}/2\mathbb{Z}}\mathbb{Z}/6\mathbb{Z}$, has no proper elementary subgroup.

\smallskip

Last, let us mention another interesting consequence of Theorem \ref{1}: if an endomorphism $\varphi$ of a virtually free group $G$ is $\forall\exists$-elementary, then $\varphi$ is an automorphism. Indeed, note that $\varphi(G)$ is a $\forall\exists$-elementary subgroup of $G$, and let us prove that $\varphi(G)=G$. Assume towards a contradiction that $\varphi(G)$ is a proper subgroup of $G$. It follows from Theorem \ref{1} that $G$ is a multiple legal large extension of $\varphi(G)$. Hence, there exists an integer $n\geq 1$ such that the abelianizations of $G$ and $\varphi(G)$ satisfy $G^{\mathrm{ab}}=\varphi(G)^{\mathrm{ab}}\times\mathbb{Z}^n$. But $\varphi(G)$ is isomorphic to $G$ since $\varphi$ is injective, hence $G^{\mathrm{ab}}\simeq G^{\mathrm{ab}}\times\mathbb{Z}^n$, which contradicts the fact that $n$ is non-zero and $G$ is finitely generated. Thus, one has $\varphi(G)=G$ and $\varphi$ is an automorphism. 

In fact, the same result holds for torsion-free hyperbolic groups: if $\varphi : G \rightarrow G$ is $\forall\exists$-elementary, then $G$ is a hyperbolic tower in the sense of Sela over $\varphi(G)\simeq G$ (see \cite{Per11}). By definition of a hyperbolic tower, $\varphi(G)$ is a quotient of $G$. Since torsion-free hyperbolic groups are Hopfian by \cite{SelHopf}, $\varphi(G)=G$ and $\varphi$ is an automorphism.

We shall prove the following result, which generalizes the previous observation. Recall that a group is said to be \emph{equationally noetherian} if every infinite system of equations $\Sigma$ in finitely many variables is equivalent to a finite subsystem of $\Sigma$.

\begin{te}\label{44}Let $G$ be a finitely generated group. Suppose that $G$ is equationally noetherian, or finitely presented and Hopfian. Then, every $\exists^{+}$-endomorphism of $G$ is an automorphism.\end{te}

\begin{rque}Note that $\forall\exists$-elementary morphisms are \emph{a fortiori} $\exists^{+}$-elementary. Note also that, contrary to $\forall\exists$-elementary morphisms, $\exists^{+}$-elementary morphisms are not injective in general.\end{rque}

As a consequence, by Proposition 2 in \cite{MR18}, a finitely generated group $G$ satisfying the hypotheses of Theorem \ref{44} above is (strongly) defined by types, and even by $\exists^{+}$-types, meaning that $G$ is characterized among finitely generated groups, up to isomorphism, by the set $\mathrm{tp}_{\exists^{+}}(G)$ of all $\exists^{+}$-types of tuples of elements of $G$. In particular, Theorem \ref{44} answers positively Problem 4 posed in \cite{MR18} and recovers several results proved in \cite{MR18}.

\section{Preliminaries}

\subsection{First-order formulas}\label{logic}

A \emph{first-order formula} in the language of groups is a finite formula using the following symbols: $\forall$, $\exists$, $=$, $\wedge$, $\vee$, $\Rightarrow$, $\neq$, $1$ (standing for the identity element), ${}^{-1}$ (standing for the inverse), $\cdot$ (standing for the group multiplication) and variables $x,y,g,z\ldots$ which are to be interpreted as elements of a group. A variable is \emph{free} if it is not bound by any quantifier $\forall$ or $\exists$. A \emph{sentence} is a formula without free variables. An \emph{existential formula} (or $\exists$\emph{-formula}) is a formula in which the symbol $\forall$ does not appear. An \emph{existential positive formula} (or $\exists^{+}$\emph{-formula}) is a formula in which the symbols $\forall$ and $\neq$ do not appear. A $\forall\exists$\emph{-formula} is a formula of the form $\theta(\bm{x}):\forall\bm{y}\exists\bm{z} \ \varphi(\bm{x},\bm{y},\bm{z})$ where $\varphi(\bm{x},\bm{y},\bm{z})$ is a quantifier-free formula, i.e.\ a boolean combination of equations and inequations in the variables of the tuples $\bm{x},\bm{y},\bm{z}$. An $\exists\forall\exists$\emph{-formula} is defined in a similar way.

\subsection{Properties relative to a subgroup}

Let $G$ be a finitely generated group, and let $H$ be a subgroup of $G$. 

\begin{de}\label{action of the pair}An \emph{action} of the pair $(G,H)$ on a simplicial tree $T$ is an action of $G$ on $T$ such that $H$ fixes a point of $T$. We always assume that the action is \emph{minimal}, which means that there is no proper subtree of $T$ invariant under the action of $G$. The tree $T$ (or the quotient graph of groups $T/G$, which is finite since the action is minimal) is called a \emph{splitting of} $(G,H)$, or a \emph{splitting of }$G$\emph{ relative to }$H$. The action is said to be \emph{trivial} if $G$ fixes a point of $T$.\end{de}

\begin{de}\label{one-ended relative to}We say that $G$ is \emph{one-ended relative to} $H$ if $G$ does not split as an amalgamated product $A\ast_C B$ or as an HNN extension $A\ast_C$ such that $C$ is finite and $H$ is contained in a conjugate of $A$ or $B$. In other words, $G$ is one-ended relative to $H$ if any action of the pair $(G,H)$ on a simplicial tree with finite edge stabilizers is trivial.\end{de}

\begin{de}\label{co-Hopfian relative to}The group $G$ is said to be \emph{co-Hopfian relative to} $H$ if every monomorphism $\varphi : G \hookrightarrow G$ that coincides with the identity on $H$ is an automorphism of $G$.\end{de}

The following result was first proved by Sela in \cite{Sel97} for torsion-free one-ended hyperbolic groups, with $H$ trivial.

\begin{te}[see \cite{And18b} Theorem 2.31]\label{coHopf}Let $G$ be a hyperbolic group, let $H$ be a subgroup of $G$. Assume that $G$ is one-ended relative to $H$. Then $G$ is co-Hopfian relative to $H$.\end{te}

\begin{rque}In \cite{And18b}, this result is stated and proved under the assumption that $H$ is finitely generated. However, Lemma \ref{perinlemme} below shows that this hypothesis is not necessary.\end{rque}

\subsection{Relative Stallings splittings}\label{chap2refchap3}

Let $G$ be a finitely generated group. Under the hypothesis that there exists a constant $C$ such that every finite subgroup of $G$ has order less than $C$, Linnell proved in \cite{Lin83} that $G$ splits as a finite graph of groups with finite edge groups and all of whose vertex groups are finite or one-ended. The group $G$ is virtually free if and only if all vertex groups are finite. Given a subgroup $H$ of $G$, Linnell's result can be generalized as follows: the pair $(G,H)$ splits as a finite graph of groups with finite edge groups such that each vertex group is finite or one-ended relative to a conjugate of $H$. Such a splitting is called a \emph{Stallings splitting of} $G$ \emph{relative to} $H$. Note that if $H$ is infinite, there exists a unique vertex group containing $H$, called the \emph{one-ended factor of} $G$ \emph{relative to} $H$. In particular, if $G$ is infinite hyperbolic and $H$ is $\exists$-elementary, then $H$ is infinite, because it satisfies the sentence $\exists x \ (x^{K_G!}\neq 1)$ where $K_G$ denotes the maximal order of an element of $G$ of finite order. As a consequence, in the context of Theorem \ref{1}, the one-ended factor of $G$ relative to $H$ is well-defined.

\subsection{The JSJ decomposition and the modular group}\label{25}

Let us denote by $\mathcal{Z}$ the class of groups that are either finite or virtually cyclic with infinite center. Let $G$ be a hyperbolic group, and let $H$ be a subgroup of $G$. Suppose that $G$ is one-ended relative to $H$. In \cite{GL16}, Guirardel and Levitt construct a splitting of $G$ relative to $H$ called the canonical JSJ splitting of $G$ over $\mathcal{Z}$ relative to $H$. In what follows, we refer to this decomposition as the $\mathcal{Z}$\emph{-JSJ splitting of} $G$ \emph{relative to} $H$. This tree $T$ enjoys particularly nice properties and is a powerful tool for studying the pair $(G,H)$. Before giving a description of $T$, let us recall briefly some basic facts about hyperbolic 2-dimensional orbifolds.

A compact connected 2-dimensional orbifold with boundary $\mathcal{O}$ is said to be \emph{hyperbolic} if it is equipped with a hyperbolic metric with totally geodesic boundary. It is the quotient of a closed convex subset $C\subset \mathbb{H}^2$ by a proper discontinuous group of isometries $G_{\mathcal{O}}\subset \mathrm{Isom}(\mathbb{H}^2)$. We denote by $p : C \rightarrow \mathcal{O}$ the quotient map. By definition, the orbifold fundamental group $\pi_1(\mathcal{O})$ of $\mathcal{O}$ is $G_{\mathcal{O}}$. We may also view $\mathcal{O}$ as the quotient of a compact orientable hyperbolic surface with geodesic boundary by a finite group of isometries. A point of $\mathcal{O}$ is \emph{singular} if its preimages in $C$ have non-trivial stabilizer. A \emph{mirror} is the image by $p$ of a component of the fixed point set of an orientation-reversing element of $G_{\mathcal{O}}$ in $C$. Singular points not contained in mirrors are \emph{conical points}; the stabilizer of the preimage in $\mathbb{H}^2$ of a conical point is a finite cyclic group consisting of orientation-preserving maps (rotations). The orbifold $\mathcal{O}$ is said to be \emph{conical} if it has no mirror.

\begin{de}\label{FBO222}A group $G$ is called a \emph{finite-by-orbifold group} if it is an extension \[1\rightarrow F\rightarrow G \rightarrow \pi_1(\mathcal{O})\rightarrow 1\]where $\mathcal{O}$ is a compact connected hyperbolic conical 2-orbifold, possibly with totally geodesic boundary, and $F$ is an arbitrary finite group called the \emph{fiber}. We call an \emph{extended boundary subgroup} of $G$ the preimage in $G$ of a boundary subgroup of the orbifold fundamental group $\pi_1(\mathcal{O})$ (for an indifferent choice of regular base point). We define in the same way \emph{extended conical subgroups}.\end{de}

\begin{de}\label{QH222}A vertex $v$ of a graph of groups is said to be \emph{quadratically hanging} (denoted by \emph{QH}) if its stabilizer $G_v$ is a finite-by-orbifold group $1\rightarrow F\rightarrow G \rightarrow \pi_1(\mathcal{O})\rightarrow 1$ such that $\mathcal{O}$ has non-empty boundary, and such that any incident edge group is finite or contained in an extended boundary subgroup of $G$. We also say that $G_v$ is QH. 
\end{de}

\begin{de}Let $G$ be a hyperbolic group, and let $H$ be a finitely generated subgroup of $G$. Let $T$ be the $\mathcal{Z}$-JSJ decomposition of $G$ relative to $H$. A vertex group $G_v$ of $T$ is said to be \emph{rigid} if it is elliptic in every splitting of $G$ over $\mathcal{Z}$ relative to $H$.\end{de}

The following proposition is crucial (see Section 6 of \cite{GL16}, Theorem 6.5 and the paragraph below Remark 9.29). We keep the same notations as in the previous definition.

\begin{prop}\label{flexible implique rigideNash}If $G_v$ is not rigid, i.e.\ if it fails to be elliptic in some splitting of $G$ over $\mathcal{Z}$ relative to $H$, then $G_v$ is quadratically hanging.\end{prop}

Proposition \ref{JSJpropNash} below summarizes the properties of the $\mathcal{Z}$-JSJ splitting relative to $H$ that are useful in the proof of Theorem \ref{1}.

\begin{prop}\label{JSJpropNash}Let $G$ be a hyperbolic group, and let $H$ be a subgroup of $G$. Suppose that $G$ is one-ended relative to $H$. Let $T$ be its $\mathcal{Z}$-JSJ decomposition relative to $H$.
\begin{itemize}
\item[$\bullet$]The tree $T$ is bipartite: every edge joins a vertex carrying a maximal virtually cyclic group to a vertex carrying a non virtually cyclic group.
\item[$\bullet$]The action of $G$ on $T$ is acylindrical in the following strong sense: if an element $g\in G$ of infinite order fixes a segment of length $\geq 2$ in $T$, then this segment has length exactly 2 and its midpoint has virtually cyclic stabilizer. 
\item Let $v$ be a vertex of $T$, and let $e,e'$ be two distinct edges incident to $v$. If $G_v$ is not virtually cyclic, then the group $\langle G_e,G_{e'}\rangle$ is not virtually cyclic. 
\item[$\bullet$]If $v$ is a QH vertex of $T$, every edge group $G_e$ of an edge $e$ incident to $v$ coincides with an extended boundary subgroup of $G_v$. Moreover, given any extended boundary subgroup $B$ of $G_v$, there exists a unique incident edge $e$ such that $G_e=B$.
\item[$\bullet$]The subgroup $H$ is contained in a rigid vertex group.
\end{itemize}
\end{prop}

\begin{rque}The rigid vertex group containing $H$ may be QH.\end{rque}

\begin{de}\label{modul}Let $G$ be a hyperbolic group and let $H$ be a subgroup of $G$. Suppose that $G$ is one-ended relative to $H$. We denote by $\mathrm{Aut}_H(G)$ the subgroup of $\mathrm{Aut}(G)$ consisting of all automorphisms whose restriction to $H$ is the conjugacy by an element of $G$. The \emph{modular group} $\mathrm{Mod}_H(G)$ of $G$ relative to $H$ is the subgroup of $\mathrm{Aut}_H(G)$ consisting of all automorphisms $\sigma$ satisfying the following conditions:
\begin{itemize}
\item[$\bullet$]the restriction of $\sigma$ to each rigid or virtually cyclic vertex group of the $\mathcal{Z}$-JSJ splitting of $G$ relative to $H$ coincides with the conjugacy by an element of $G$,
\item[$\bullet$]the restriction of $\sigma$ to each finite subgroup of $G$ coincides with the conjugacy by an element of $G$,
\item[$\bullet$]$\sigma$ acts trivially on the underlying graph of the $\mathcal{Z}$-JSJ splitting relative to $H$.
\end{itemize}
\end{de}

We will need the following result.


\begin{te}\label{short}Let $G$ be a hyperbolic group, let $H$ be a subgroup of $G$ and let $U$ be the one-ended factor of $G$ relative to $H$. There exist a finite subset $F\subset U\setminus \lbrace 1\rbrace$ and a finitely generated subgroup $H'\subset H$ such that, for every non-injective homomorphism $\varphi : U \rightarrow G$ that coincides with the identity on $H'$ up to conjugation, there exists an automorphism $\sigma\in\mathrm{Mod}_H(U)$ such that $\ker(\varphi\circ\sigma)\cap F\neq \varnothing$.\end{te}

\begin{proof}This result is stated and proved in \cite{And18b} under the assumption that $H$ is finitely generated (see Theorem 2.32), in which case one can take $H'=H$. We only give a brief sketch of how the proof can be adapted if $H$ is not assumed to be finitely generated. In \cite{And18b}, the assumption that $H$ is finitely generated is only used in the proof of  Proposition 2.27 in order to ensure that the group $H$ fixes a point in a certain real tree $T$ with virtually cyclic arc stabilizers (namely the tree obtained by rescaling the metric of a Cayley graph of $G$ by a given sequence of positive real numbers going to infinity). Let $\lbrace h_1,h_2,\ldots\rbrace$ be a generating set for $H$, and let $H_n$ be the subgroup of $H$ generated by $\lbrace h_1,\ldots,h_n\rbrace$. If $H$ is not finitely generated, then there exists an integer $n_0$ such that, for all $n\geq n_0$, the subgroup $H_n$ is not virtually cyclic. It follows that all $H_n$ fix the same point of $T$ for $n\geq n_0$, which proves that $H$ is elliptic in $T$. Hence, one can just take $H'=H_{n_0}$.\end{proof}

\subsection{Related homomorphisms and preretractions}\label{322}

We denote by $\mathrm{ad}(g)$ the inner automorphism $h\mapsto ghg^{-1}$.

\begin{de}[Related homomorphisms]\label{reliés222}
Let $G$ be a hyperbolic group and let $H$ be a subgroup of $G$. Assume that $G$ is one-ended relative to $H$. Let $G'$ be a group. Let $\Lambda$ be the $\mathcal{Z}$-JSJ splitting of $G$ relative to $H$. Let $\varphi$ and $\varphi'$ be two homomorphisms from $G$ to $G'$. We say that $\varphi$ and $\varphi'$ are $\mathcal{Z}$-\emph{JSJ-related} or \emph{$\Lambda$-related} if the following two conditions hold:
\begin{itemize}
\item[$\bullet$]for every vertex $v$ of $\Lambda$ such that $G_v$ is rigid or virtually cyclic, there exists an element $g_v\in G'$ such that \[{\varphi'}_{\vert G_v}=\mathrm{ad}({g_v})\circ \varphi_{\vert G_v};\]
\item[$\bullet$]for every finite subgroup $F$ of $G$, there exists an element $g\in G'$ such that \[{\varphi'}_{\vert F}=\mathrm{ad}({g})\circ \varphi_{\vert F}.\]
\end{itemize}
\end{de}

\begin{de}[Preretraction]\label{pre222}
Let $G$ be a hyperbolic group, and let $H$ be a subgroup of $G$. Assume that $G$ is one-ended relative to $H$. Let $\Lambda$ be the $\mathcal{Z}$-JSJ splitting of $G$ relative to $H$. A $\mathcal{Z}$-\emph{JSJ-preretraction} or $\Lambda$\emph{-preretraction} of $G$ is an endomorphism of $G$ that is $\Lambda$-related to the identity map. More generally, if $G$ is a subgroup of a group $G'$, a preretraction from $G$ to $G'$ is a homomorphism $\Lambda$-related to the inclusion of $G$ into $G'$. Note that a $\Lambda$-preretraction coincides with a conjugacy on $H$, since $H$ is contained in a rigid vertex group of $\Lambda$.
\end{de}

The following easy lemma shows that being $\Lambda$-related can be expressed in first-order logic. This lemma is stated and proved in \cite{And18b} (see Lemma 2.22) under the assumption that $H$ is finitely generated, but this hypothesis is not used in the proof.

\begin{lemme}\label{deltarelies222}
Let $G$ be a hyperbolic group and let $H$ be a subgroup of $G$. Assume that $G$ is one-ended relative to $H$. Let $G'$ be a group. Let $\Lambda$ be the $\mathcal{Z}$-JSJ splitting of $G$ relative to $H$. Let $\lbrace g_1,\ldots ,g_n\rbrace$ be a generating set of $G$. There exists an existential formula $\theta(x_1,\ldots , x_{2n})$ with $2n$ free variables such that, for every $\varphi,\varphi'\in \mathrm{Hom}(G,G')$, $\varphi$ and $\varphi'$ are $\Lambda$-related if and only if $G'$ satisfies $\theta (\varphi(g_1),\ldots , \varphi(g_n),\varphi'(g_1),\ldots ,\varphi'(g_n))$.
\end{lemme}

The proof of the following lemma is identical to that of Proposition 7.2 in \cite{And18}.

\begin{lemme}\label{lemmeperin2}Let $G$ be a hyperbolic group. Suppose that $G$ is one-ended relative to a subgroup $H$. Let $\Lambda$ be the $\mathcal{Z}$-JSJ splitting of $G$ relative to $H$. Let $\varphi$ be a $\Lambda$-preretraction of $G$. If $\varphi$ sends every QH vertex group of $\Lambda$ isomorphically to a conjugate of itself, then $\varphi$ is injective.
\end{lemme}

\subsection{Centered graph of groups}\label{centered}

\begin{de}[Centered graph of groups]\label{graphecentre222}A graph of groups over $\mathcal{Z}$, with at least two vertices, is said to be \emph{centered} if the following conditions hold:
\begin{itemize}
\item[$\bullet$]the underlying graph is bipartite, with a particular QH vertex $v$ such that every vertex different from $v$ is adjacent to $v$;
\item[$\bullet$]every stabilizer $G_e$ of an edge incident to $v$ coincides with an extended boundary subgroup or with an extended conical subgroup of $G_v$ (see Definition \ref{FBO222});
\item[$\bullet$]given any extended boundary subgroup $B$, there exists a unique edge $e$ incident to $v$ such that $G_e$ is conjugate to $B$ in $G_v$;
\item[$\bullet$]if an element of infinite order fixes a segment of length $\geq 2$ in the Bass-Serre tree of the splitting, then this segment has length exactly 2 and its endpoints are translates of $v$.
\end{itemize}
The vertex $v$ is called \emph{the central vertex}.
\end{de}

\begin{figure}[!h]
\includegraphics[scale=0.4]{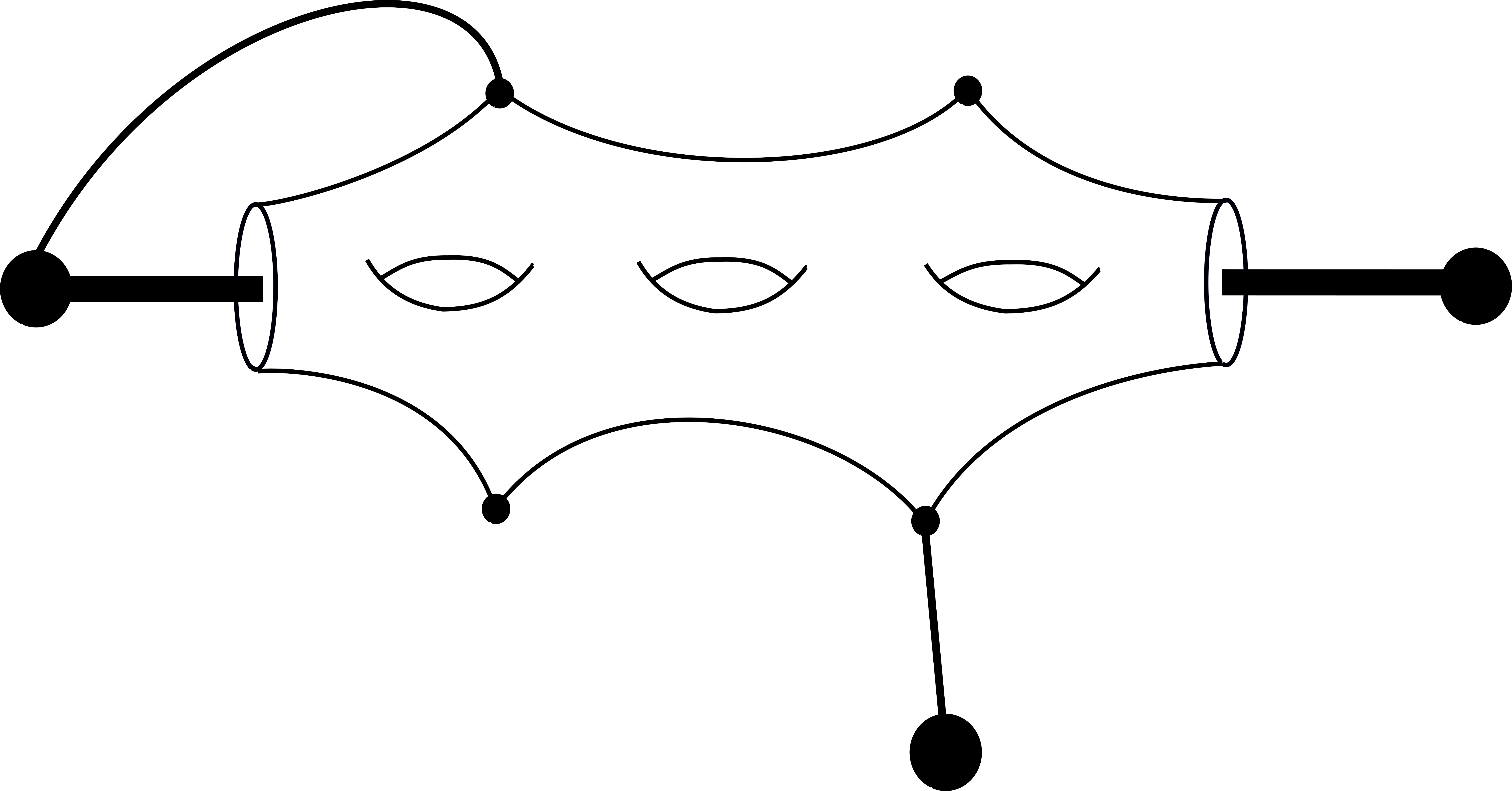}
\caption{A centered graph of groups. Edges with infinite stabilizer are depicted in bold.}
\end{figure}


We also need to define relatedness and preretractions in the context of centered graphs of groups.

\begin{de}[Related homomorphisms]\label{reliés2222}
\normalfont
Let $G$ and $G'$ be two groups. Let $H$ be a subgroup of $G$. Suppose that $G$ has a centered splitting $\Delta$, with central vertex $v$. Suppose that $H$ is contained in a non-central vertex of $\Delta$. Let $\varphi$ and $\varphi'$ be two homomorphisms from $G$ to $G'$. We say that $\varphi$ and $\varphi'$ are \emph{$\Delta$-related} (relative to $H$) if the following two conditions hold:
\begin{itemize}
\item[$\bullet$]for every vertex $w\neq v$, there exists an element $g_w\in G'$ such that \[{\varphi'}_{\vert G_w}=\mathrm{ad}(g_w)\circ \varphi_{\vert G_w};\]
\item[$\bullet$]for every finite subgroup $F$ of $G$, there exists an element $g\in G'$ such that \[{\varphi'}_{\vert F}=\mathrm{ad}(g)\circ \varphi_{\vert F}.\]
\end{itemize}
\end{de}

\begin{de}[Preretraction]\label{special}Let $G$ be a hyperbolic group, let $H$ be a subgroup of $G$, and let $\Delta$ be a centered splitting of $G$. Let $v$ be the central vertex of $\Delta$. Suppose that $H$ is contained in a non-central vertex of $\Delta$. An endomorphism $\varphi$ of $G$ is called a \emph{$\Delta$-preretraction} (relative to $H$) if it is $\Delta$-related to the identity of $G$ in the sense of the previous definition. A $\Delta$-preretraction is said to be \emph{non-degenerate} if it does not send $G_v$ isomorphically to a conjugate of itself.\end{de}

\section{Elementary subgroups of virtually free groups}

In this section, we prove Theorem \ref{1}. Recall that this theorem claims that if $G$ is a virtually free group and $H$ is a proper $\forall\exists$-elementary subgroup of $G$, then $G$ is a multiple legal large extension of $H$.

\subsection{Elementary subgroups are one-ended factors}

As a first step, we will prove the following result.

\begin{prop}\label{256}Let $G$ be a virtually free group. Let $H$ be a subgroup of $G$. If $H$ is $\forall\exists$-elementary, then $H$ coincides with the one-ended factor of $G$ relative to $H$. In other words, $H$ appears as a vertex group in a splitting of $G$ over finite groups.\end{prop}

The proof of Proposition \ref{256}, which is inspired from \cite{Per11}, consists in showing that if $G$ is a hyperbolic group and $H$ is strictly contained in the one-ended factor of $G$ relative to $H$, then there exists a centered splitting $\Delta$ of $G$ relative to $H$, and a non-degenerate $\Delta$-preretraction of $G$ (see Lemmas \ref{2} and \ref{lemmelemme} below). However, if $G$ is virtually free, Lemma \ref{cyclic2} below shows that such a preretraction cannot exist.

We shall prove Proposition \ref{256} after establishing a series of preliminary lemmas. The following result is a generalization of Lemma 4.20 in \cite{Per11}. Recall that all group actions on trees considered in this paper are assumed to be minimal (see Definition \ref{action of the pair}). As a consequence, trees have no vertex of valence 1. We say that a tree $T$ endowed with an action of a group $G$ is \emph{non-redundant} if there exists no valence 2 vertex $v$ such that both boundary monomorphisms into the vertex group $G_v$ are isomorphisms.

\begin{lemme}\label{perinlemme}Let $G$ be a finitely generated group, and let $H$ be a subgroup of $G$. Suppose that $G$ is one-ended relative to $H$ and that there is a constant $C$ such that every finite subgroup of $G$ has order at most $C$. Then there exists a finitely generated subgroup $H''$ of $H$ such that $G$ is one-ended relative to $H''$.
\end{lemme}

\begin{proof}
Let $\lbrace h_1,h_2,\ldots \rbrace$ be a generating set for $H$, possibly infinite. For every integer $n\geq 1$, let $H_n$ be the subgroup of $H$ generated by $\lbrace h_1,\ldots,h_n \rbrace$. By Theorem 1 in \cite{Wei12}, there is a maximum number $m_n$ of orbits of edges in a non-redundant splitting of $G$ relative to $H_n$ over finite groups. Let $T_n$ be such a splitting with $m_n$ orbits of edges, and let $G_n$ be the vertex group of $T_n$ containing $H_n$. 


We shall prove that $G_{n+1}$ is contained in $G_n$ for all $n$ sufficiently large. First, note that the sequence of integers $(m_n)_{n\in\mathbb{N}}$ is non-increasing, because $T_{n+1}$ is a splitting of $G$ relative to $H_n$. In particular, there exists an integer $n_0$ such that $m_n=m_{n+1}$ for every $n\geq n_0$. We claim that $G_{n+1}$ is elliptic in $T_n$. Otherwise, there exists a non-trivial splitting $G_{n+1}=A\ast_CB$ or $G_{n+1}=A\ast_C$ with $C$ finite and $H_n\subset A$, and one gets a non-redundant splitting of $G$ relative to $H_n$ over finite groups with $m_{n+1}+1=m_n+1$ edges by replacing the vertex group $G_{n+1}$ of the graph of groups $T_{n+1}/G$ with the previous one-edge splitting of $G_{n+1}$, which contradicts the definition of $m_n$.

Hence, for $n\geq n_0$, one has $G_{n}\subset G_{n_0}$. In particular, $G_{n_0}$ contains $H_n$ for every integer $n$. Thus, $G_{n_0}$ contains $H$. Since $G$ is assumed to be one-ended relative to $H$, one has $G=G_{n_0}$ and one can take $H''=H_{n_0}$.
\end{proof}

We will need the following well-known result in the proof of Lemma \ref{2} below.

\begin{prop}[\cite{Bow98}, Proposition 1.2]\label{bowNash}If a hyperbolic group splits over quasi-convex subgroups, then every vertex group is quasi-convex (hence hyperbolic).
\end{prop}

\begin{lemme}\label{2}Let $G$ be a hyperbolic group. Let $H$ be a $\forall\exists$-elementary subgroup of $G$. Let $U$ be the one-ended factor of $G$ relative to $H$. Let $\Lambda$ be the $\mathcal{Z}$-JSJ splitting of $U$ relative to $H$. If $H$ is strictly contained in $U$, then there exists a non-injective $\Lambda$-preretraction $U\rightarrow G$.\end{lemme}

\begin{proof}Let $H'$ be the finitely generated subgroup of $H$ given by Theorem \ref{short} and let $H''$ be the finitely generated subgroup of $H$ given by Lemma \ref{perinlemme} above. Let $H_0$ be the finitely generated subgroup of $H$ generated by $H'\cup H''$.

Let us prove that every morphism $\varphi:U\rightarrow H$ whose restriction to $H_0$ coincides with the identity is non-injective. First, note that $U$ is one-ended relative to $H_0$ (since it is one-ended related to $H''$ which is contained in $H_0$), and that $U$ is hyperbolic by Proposition \ref{bowNash} above. Therefore, by Theorem \ref{coHopf}, $U$ is co-Hopfian relative to $H_0$. Hence, a putative monomorphism $\varphi:U\rightarrow H\subset U$ whose restriction to $H_0$ coincides with the identity is surjective, viewed as an endomorphism of $U$. But $\varphi(U)$ is contained in $H$, which shows that $U=\varphi(U)$ is contained in $H$. It's a contradiction since $H$ is stricly contained in $U$, by assumption. 

We proved in the previous paragraph that every morphism $\varphi:U\rightarrow H$ whose restriction to $H_0$ coincides with the identity is non-injective. Therefore, by Theorem \ref{short}, for every morphism $\varphi:U\rightarrow H$ whose restriction to $H_0$ (which contains $H'$) coincides with the identity, there exists an automorphism $\sigma\in\mathrm{Mod}_{H}(U)$ such that $\varphi\circ \sigma$ kills an element of the finite set $F\subset U\setminus \lbrace 1\rbrace$ given by Theorem \ref{short}. In addition, note that the morphisms $\varphi\circ \sigma$ and $\varphi$ are $\Lambda$-related (see Definition \ref{reliés2222}). Hence, for every morphism $\varphi:U\rightarrow H$ whose restriction to $H_0$ coincides with the identity, there exists a morphism $\varphi':U\rightarrow H$ that kills an element of the finite set $F$, and which is $\Lambda$-related to $\varphi$. We will see that this statement $(\star)$ is expressible by means of a $\forall\exists$-sentence with constants in $H$.

Let $U=\langle u_1,\ldots,u_n \ \vert \ R(u_1,\ldots,u_n)=1\rangle$ be a finite presentation of $U$. Let $\lbrace h_1,\ldots,h_p\rbrace$ be a finite generating set for $H_0$. For every $1\leq i\leq p$, the element $h_i$ can be written as a word $w_i(u_1,\ldots,u_n)$. Likewise, one can write $F=\lbrace v_1(u_1,\ldots,u_n),\ldots,v_k(u_1,\ldots,u_n)\rbrace$.

Observe that there is a one-to-one correspondence between the set of homomorphisms $\mathrm{Hom}(U,H)$ and the set of solutions in $H^n$ of the system of equations $R(x_1,\ldots,x_n)=1$. The group $H$ satisfies the following $\forall\exists$-formula, expressing the statement $(\star)$:
\begin{small}
\begin{align*}
&\mu(h_1,\ldots,h_p): \ \forall x_1\ldots\forall x_n \ \left(R(x_1,\ldots,x_n)=1 \ \wedge \ \bigwedge_{i=1}^p w_i(x_1,\ldots,x_n)=h_i\right)\\&\Rightarrow \left(\exists x'_1\ldots\exists x'_n \ R(x'_1,\ldots,x'_n)=1 \ \wedge \ \theta(x_1,\ldots,x_n,x'_1,\ldots,x'_n)=1 \wedge \bigvee_{i=1}^k v_i(x'_1,\ldots,x'_n)= 1 \right)
\end{align*}
\end{small}
where $\theta$ is the formula given by Lemma \ref{deltarelies222}, expressing that the homomorphisms $\varphi$ and $\varphi'$ defined by $h_i\mapsto x_i$ and $\varphi' : h_i \mapsto x'_i$ are $\Lambda$-related, where $\Lambda$ denotes the $\mathcal{Z}$-JSJ splitting of $U$ relative to $H$.

Since $H$ is $\forall\exists$-elementary (as a subgroup of $G$), the group $G$ satisfies $\mu(h_1,\ldots,h_p)$ as well. For $x_i=u_i$ for $1\leq i\leq p$, the interpretation of $\mu(h_1,\ldots,h_p)$ in $G$ provides a tuple $(g_1,\ldots,g_n)\in G^n$ such that the application $p : U \rightarrow G$ defined by $u_i \mapsto g_i$ for every $1\leq i\leq p$ is a homomorphism, is $\Lambda$-related to the inclusion of $U$ into $G$ (see Definition \ref{reliés2222}), and kills an element of $F$. As a conclusion, $p$ is a non-injective $\Lambda$-preretraction from $U$ to $G$ (see Definition \ref{pre222}).\end{proof}

The following easy lemma is proved in \cite{And18b} (see Lemma 4.5).

\begin{lemme}\label{petitlemme}
Let $G$ be a group endowed with a splitting over finite groups. Let $T_G$ denote the Bass-Serre tree associated with this splitting. Let $U$ be a group endowed with a splitting over infinite groups, and let $T_U$ be the associated Bass-Serre tree. If $p:U\rightarrow G$ is a homomorphism injective on edge groups of $T_U$, and such that $p(U_v)$ is elliptic in $T_G$ for every vertex $v$ of $T_U$, then $p(U)$ is elliptic in $T_G$.
\end{lemme}

\begin{lemme}\label{lemmelemme}Let $G$ be a hyperbolic group. Let $H$ be a subgroup of $G$. Let $U$ be the one-ended factor of $G$ relative to $H$. Let $\Lambda$ be the $\mathcal{Z}$-JSJ splitting of $G$ relative to $H$. Suppose that there exists a non-injective $\Lambda$-preretraction $p:U\rightarrow G$. Then there exists a centered splitting of $G$ relative to $H$, called $\Delta$, and a non-degenerate $\Delta$-preretraction of $G$.
\end{lemme}

\begin{proof}First, we will prove that there exists a QH vertex $x$ of $\Lambda$ such that $U_x$ is not sent isomorphically to a conjugate of itself by $p$. Assume towards a contradiction that this claim is false, i.e.\ that each stabilizer $U_x$ of a QH vertex $x$ of $\Lambda$ is sent isomorphically to a conjugate of itself by $p$. We claim that $p(U)$ is contained in a conjugate of $U$. Let $\Gamma$ be a Stallings splitting of $G$ relative to $H$. By definition of $U$, there exists a vertex $u$ of the Bass-Serre tree $T$ of $\Gamma$ such that $G_u=U$. First, let us check that the hypotheses of Lemma \ref{petitlemme} are satisfied:
\begin{enumerate}
\item by definition, $\Gamma$ is a splitting of $G$ over finite groups, and $\Lambda$ is a splitting of $U$ over infinite groups;
\item $p$ is injective on edge groups of $\Lambda$ (as a $\Lambda$-preretraction);
\item if $x$ is a QH vertex of $\Lambda$, then $p(U_x)$ is conjugate to $U_x$ by assumption. In particular, $p(U_x)$ is contained in a conjugate of $U$ in $G$. As a consequence, $p(U_x)$ is elliptic in $T$ (more precisely, it fixes a translate of the vertex $u$ of $T$ such that $G_u=U$). If $x$ is a non-QH vertex of $\Lambda$, then $p(U_x)$ is conjugate to $U_x$ by definition of a $\Lambda$-preretraction. In particular, $p(U_x)$ is elliptic in $T$.
\end{enumerate}
By Lemma \ref{petitlemme}, $p(U)$ is elliptic in $T$. It remains to prove that $p(U)$ is contained in a conjugate of $U$. Observe that $U$ is not finite-by-(closed orbifold), as a virtually free group. Therefore, there exists at least one non-QH vertex $x$ in $\Lambda$. Moreover, since $p$ is inner on non-QH vertices of $\Lambda$, there exists an element $g\in G$ such that $p(U_x)=gU_xg^{-1}$. Hence, $p(U)\cap gUg^{-1}$ is infinite, which proves that $p(U)$ is contained in $gUg^{-1}$ since edge groups of the Bass-Serre tree $T$ of $\Gamma$ are finite.

Now, up to composing $p$ by the conjugation by $g^{-1}$, one can assume that $p$ is an endomorphism of $U$. By Lemma \ref{lemmeperin2}, $p$ is injective. This is a contradiction since $p$ is non-injective by hypothesis. Hence, we have proved that there exists a QH vertex $x$ of $\Lambda$ such that $U_x$ is not sent isomorphically to a conjugate of itself by $p$.

Then, we refine $\Gamma$ by replacing the vertex $u$ fixed by $U$ by the $\mathcal{Z}$-JSJ splitting $\Lambda$ of $U$ relative to $H$ (which is possible since edge groups of $\Gamma$ adjacent to $u$ are finite, ans thus are elliptic in $\Lambda$). With a little abuse of notation, we still denote by $x$ the vertex of $\Gamma$ corresponding to the QH vertex $x$ of $\Lambda$. Then, we collapse to a point every connected component of the complement of $\mathrm{star}(x)$ in $\Gamma$ (where $\mathrm{star}(x)$ stands for the subgraph of $\Gamma$ constituted of $x$ and all its incident edges). The resulting graph of groups, denoted by $\Delta$, is non-trivial. One easily sees that $\Delta$ is a centered splitting of $G$, with central vertex $x$.

The homomorphism $p:U\rightarrow G$ is well-defined on $G_x$ because $G_x=U_x$ is contained in $U$. Moreover, $p$ restricts to a conjugation on each stabilizer of an edge $e$ of $\Delta$ incident to $x$. Indeed, either $e$ is an edge coming from $\Lambda$, either $G_e$ is a finite subgroup of $U$; in each case, $p_{\vert G_e}$ is a conjugation since $p$ is $\Lambda$-related to the inclusion of $U$ into $G$. Now, one can define an endomorphism $\varphi : G \rightarrow G$ that coincides with $p$ on $G_x=U_x$ and coincides with a conjugation on every vertex group $G_y$ of $\Gamma$, with $y\neq x$. By induction on the number of edges of $\Gamma$, it is enough to define $\varphi$ in the case where $\Gamma$ has only one edge. If $G=U_x\ast_C B$ with $p_{\vert C}=\mathrm{ad}(g)$, one defines $\varphi : G \rightarrow G$ by $\varphi_{\vert U_x}=p$ and $\varphi_{\vert B}=\mathrm{ad}(g)$. If $G=U_x\ast_C=\langle U_x,t \ \vert \ tct^{-1}=\alpha(c), \forall c\in C\rangle$ with $p_{\vert C}=\mathrm{ad}(g_1)$ and $p_{\vert \alpha{(C)}}=\mathrm{ad}(g_2)$, one defines $\varphi : G \rightarrow G$ by $\varphi_{\vert U_x}=p$ and $\varphi(t)=g_2^{-1}tg_1$.

The endomorphism $\varphi$ defined above is $\Delta$-related to the identity of $G$ (in the sense of Definition \ref{reliés2222}), and $\varphi$ does not send $G_x$ isomorphically to a conjugate of itself. Hence, $\varphi$ is a non-degenerate $\Delta$-preretraction of $G$ (see Definition \ref{special}).\end{proof}

The following result is proved in \cite{And18b} (Lemma 4.4).

\begin{lemme}\label{cyclic2}Let $G$ be a virtually free group, and let $\Delta$ be a centered splitting of $G$. Then $G$ has no non-degenerate $\Delta$-preretraction.\end{lemme}

We can now prove Proposition \ref{256}.

\medskip

\begin{proofPROP}Let $U$ be the one-ended factor of $G$ relative to $H$. Assume towards a contradiction that $H$ is strictly contained in $U$. Then by Proposition \ref{2}, there exists a non-injective preretraction $U\rightarrow G$ (with respect to the $\mathcal{Z}$-JSJ splitting of $U$ relative to $H$). By Lemma \ref{lemmelemme}, there exists a centered splitting $\Delta$ of $G$ relative to $H$ such that $G$ has a non-degenerate $\Delta$-preretraction. This contradicts Lemma \ref{cyclic2}. Hence, $H$ is equal to $U$.\end{proofPROP}


\subsection{Proof of Theorem \ref{1}}

Recall that Theorem \ref{1} claims that if $H$ is a $\forall\exists$-elementary proper subgroup of a virtually free group $G$, then $G$ is a multiple legal large extension of $H$. Before proving this result, we will define five numbers associated with a hyperbolic group, which are encoded into its $\forall\exists$-theory (see Lemma \ref{5inv} below).

\begin{de}\label{5nombres}Let $G$ be a hyperbolic group. We associate to $G$ the following five integers:
\begin{itemize}
\item[$\bullet$]the number $n_1(G)$ of conjugacy classes of finite subgroups of $G$,
\item[$\bullet$]the sum $n_2(G)$ of $\vert \mathrm{Aut}_G(C_k)\vert$ for $1\leq k\leq n_1(G)$, where the $C_k$ are representatives of the conjugacy classes of finite subgroups of $G$, and \[\mathrm{Aut}_G(C_k)=\lbrace \alpha\in \mathrm{Aut}(C_k) \ \vert \ \exists g\in N_G(C_k), \ \mathrm{ad}(g)_{\vert C}=\alpha\rbrace,\]
\item[$\bullet$]the number $n_3(G)$ of conjugacy classes of finite subgroups $C$ of $G$ such that $N_G(C)$ is infinite virtually cyclic,
\item[$\bullet$]the number $n_4(G)$ of conjugacy classes of finite subgroups $C$ of $G$ such that $N_G(C)$ is not virtually cyclic (finite or infinite),
\item[$\bullet$]the number $n_5(G)$ of conjugacy classes of finite subgroups $C$ of $G$ such that $N_G(C)$ is not virtually cyclic (finite or infinite) and $E_G(N_G(C))\neq C$.
\end{itemize}
\end{de}

The following lemma shows that these five numbers are preserved under $\forall\exists$-equivalence. Its proof is quite straightforward and is postponed after the proof of Theorem \ref{1}.

\begin{lemme}\label{5inv}Let $G$ and $G'$ be two hyperbolic groups. Suppose that $\mathrm{Th}_{\forall\exists}(G)=\mathrm{Th}_{\forall\exists}(G')$. Then $n_i(G)=n_i(G')$, for $1\leq i\leq 5$.\end{lemme}

Theorem \ref{1} will be an easy consequence of the following result.

\begin{prop}\label{concon}Let $G$ be a virtually free group. Let $H$ be a proper subgroup of $G$. Suppose that the following three conditions are satisfied:
\begin{enumerate}
\item $n_i(H)=n_i(G)$ for all $1\leq i\leq 5$,
\item $H$ appears as a vertex group in a splitting of $G$ over finite groups,
\item two finite subgroups of $H$ are conjugate in $H$ if and only if they are conjugate in $G$.
\end{enumerate}Then $G$ is a multiple legal large extension of $H$ (see Definition \ref{legal}).
\end{prop}

\begin{proof}First, note that the equality $n_4(G)=n_4(G)$ implies that $H$ is non virtually cyclic. Indeed, if $H$ is virtually cyclic, then $n_4(H)=0$, whereas $n_4(G)$ is greater than $1$ since $N_G(\lbrace 1\rbrace)=G$ is not virtually cyclic by assumption.

Let $T$ be the Bass-Serre tree of the splitting of $G$ given by the second condition. Up to refining this splitting, one can assume without loss of generality that the vertex groups of $T$ which are not conjugate to $H$ are finite. In other words, $T$ is a Stallings splitting of $G$ relative to $H$, in which $H$ is a vertex group by assumption. Moreover, up to collasping some edges, one can assume that $T$ is \emph{reduced}, which means that if $e=[v,w]$ is an edge of $T$ such that $G_e=G_v=G_w$, then $v$ and $w$ are in the same orbit. We denote by $\Gamma$ the quotient graph of groups $T/G$.

We will deduce from the third condition that the underlying graph of $\Gamma$ has only one vertex. Assume towards a contradiction that the Bass-Serre tree $T$ of $\Gamma$ has at least two orbits of vertices. Hence, there is a vertex $v$ of $T$ which is not in the orbit of the vertex $v_H$ fixed by $H$. By definition of $\Gamma$, the vertex stabilizer $G_v$ is finite. Thus, there exists an element $g\in G$ such that $gG_vg^{-1}$ is contained in $H$. Therefore, $G_v$ stabilizes the path of edges in $T$ between the vertices $v$ and $g^{-1}v_H$. It follows that $G_v$ coincides with the stabilizer of an edge incident to $v$ in $T$, which contradicts the assumption that $T$ is reduced.

Hence, the underlying graph of $\Gamma$ is a rose, and the central vertex group of $\Gamma$ is $H$. Moreover, edge stabilizers of $\Gamma$ are finite. In other words, there exist pairs of finite subgroups $(C_1,C'_1),\ldots,(C_n,C'_n)$ of $H$, together with automorphisms $\alpha_1\in\mathrm{Isom}(C_1,C'_1),\ldots,\alpha_n\in\mathrm{Isom}(C_n,C'_n)$ such that $G$ has the following presentation:\[G=\langle H,t_1,\ldots,t_n \ \vert \ \mathrm{ad}(t_i)_{\vert C_i}=\alpha_i, \ \forall i\in\llbracket 1,n\rrbracket\rangle.\]
By assumption, the integers $n_i(G)$ and $n_i(H)$ are equal, for $1\leq i\leq 5$. From the equality $n_1(G)=n_1(H)$, one deduces immediately that the finite groups $C_i$ and $C'_i$ are conjugate in $H$ for every integer $i\in \lbrace 1,\ldots ,n\rbrace$. Therefore, one can assume without loss of generality that $C'_i=C_i$. 

Note that for every finite subgroup $C$ of $H$, the group $ \mathrm{Aut}_{H}(C)$ is contained in $\mathrm{Aut}_{G}(C)$. Thus, the equality $n_2(G)=n_2(H)$ guarantees that $\mathrm{Aut}_{H}(C_i)$ is in fact equal to $\mathrm{Aut}_{G}(C_i)$, for every $1\leq i\leq n$. Hence, since the automorphism $\mathrm{ad}(t_i)_{\vert C_i}$ of $C_i$ belongs to $\mathrm{Aut}_{G}(C_i)$, there exists an element $h_i\in N_H(C_i)$ such that $\mathrm{ad}(h_i)_{\vert C_i}=\mathrm{ad}(t_i)_{\vert C_i}$. Up to replacing $t_i$ with $t_ih_i^{-1}$, the group $G$ has the following presentation:\[G=\langle H,t_1,\ldots,t_n \ \vert \ \mathrm{ad}(t_i)_{\vert C_i}=\mathrm{id}_{C_i}, \ \forall i\in\llbracket 1,n\rrbracket\rangle.\]
In order to prove that $G$ is a multiple legal large extension of $H$ (see Definition \ref{legal}), it remains to prove that the following two conditions hold, for every integer $1\leq i\leq n$:
\begin{enumerate}
\item the normalizer $N_H(C_i)$ is non virtually cyclic (finite or infinite),
\item and the finite group $E_H(N_H(C_i))$ coincides with $C_i$.
\end{enumerate}

The equalities $n_3(G)=n_3(H)$ and $n_4(G)=n_4(H)$ ensure that $N_H(C_i)$ is not virtually cyclic. Indeed, if $N_H(C_i)$ were finite, then $N_{G}(C_i)$ would be infinite virtually cyclic and $n_3(G)$ would be at least $n_3(H)+1$; similarly, if $N_H(C_i)$ were infinite virtually cyclic, then $N_{G}(C_i)$ would be non virtually cyclic and $n_4(G)\geq n_4(H)+1$. Hence, the first condition above is satisfied. 

Last, it follows from the equality $n_5(G)=n_5(H)$ that the finite group $E_H(N_H(C_i))$ coincides with $C_i$, otherwise $n_5(G)\geq n_5(H)+1$, since $E_{G}(N_{G}(C_i))=C_i$. Thus, the second condition above holds. As a conclusion, $G$ is a multiple legal large extension of $H$ in the sense of Definition \ref{legal}.
\end{proof}

We can now prove Theorem \ref{1}.

\medskip

\begin{proofTH1}Let $G$ be a virtually free group, and let $H$ be a $\forall\exists$-elementary subgroup of $G$. In particular, $G$ and $H$ have the same $\forall\exists$-theory. It follows from Lemma \ref{5inv} that $n_i(H)$ is equal to $n_i(G)$ for all $1\leq i\leq 5$. Hence, the first condition of Proposition \ref{concon} holds.

By Proposition \ref{256}, $H$ is a vertex group in a splitting of $G$ over finite groups, which means that the second condition of Proposition \ref{concon} is satisfied.

It remains to check the third condition of Proposition \ref{concon}, namely that two finite subgroups of $H$ are conjugate in $H$ if and only if they are conjugate in $G$. First, recall that $H$ and $G$ have the same number of conjugacy classes of finite subgroups, since $n_1(G)=n_1(H)$. Then, the conclusion follows from the following observation: if two finite subgroups $A=\lbrace a_1,\ldots ,a_m\rbrace$ and $B=\lbrace b_1,\ldots,b_m\rbrace$ of $H$ are not conjugate in $H$, then they are not conjugate in $G$. Indeed, $H$ satisfies the following universal formula: \[\theta(a_1,\ldots,a_m,b_1,\ldots,b_m):\forall x \ \bigvee_{i=1}^m\bigwedge_{j=1}^m xa_ix^{-1}\neq b_j.\]Since $H$ is $\forall\exists$-elementary (in particular $\forall$-elementary), $G$ satisfies this sentence as well. Therefore, $A$ and $B$ are not conjugate in $G$.\end{proofTH1}

It remains to prove Lemma \ref{5inv}. First, recall that if $G$ is hyperbolic and $g\in G$ has infinite order, there is a unique maximal virtually cyclic subgroup of $G$ containing $g$, denoted by $M(g)$. More precisely, $M(g)$ is the stabilizer of the pair of fixed points of $g$ on the boundary $\partial_{\infty} G$ of $G$. If $h$ and $g$ are two elements of infinite order, either $M(h)=M(g)$ or the intersection $M(h)\cap M(g)$ is finite; in the latter case, the subgroup $\langle h,g\rangle$ is not virtually cyclic. Let $K_G$ denote the maximum order of an element of $G$ of finite order. One can see that an element $g\in G$ has infinite order if and only if $g^{K_G!}$ is non-trivial, and that if $g$ and $h$ have infinite order, then $M(g)=M(h)$ if and only if the commutator $[g^{K!},h^{K!}]$ is trivial. In other words, the subgroup $\langle g,h\rangle$ is virtually cyclic if and only if $[g^{K!},h^{K!}]=1$.

\medskip

\begin{proofL}Let us denote by $K_G$ the maximal order of a finite subgroup of $G$. Since $G$ and $G'$ have the same existential theory, we have $K_G=K_{G'}$. Let $n\geq 1$ be an integer. If $n_1(G)\geq n$, then the following $\exists\forall$-sentence, written in natural language for convenience of the reader and denoted by $\theta_{1,n}$, is satisfied by $G$: there exist $n$ finite subgroups $C_1,\ldots ,C_n$ of $G$ such that, for every $g\in G$ and $1\leq i\neq j\leq n$, the groups $gC_ig^{-1}$ and $C_j$ are distinct. Since $G$ and $G'$ have the same $\exists\forall$-theory, the sentence $\theta_{1,n}$ is satisfied by $G'$ as well. As a consequence, $n_1(G')\geq n$. It follows that $n_1(G')\geq n_1(G)$. By symmetry, we have $n_1(G)=n_1(G')$.

\smallskip

In the rest of the proof, we give similar sentences $\theta_{2,n},\ldots,\theta_{5,n}$ such that the following series of equivalences hold: $n_i(G)\geq n \Leftrightarrow G$ satisfies $\theta_{i,n} \Leftrightarrow$ $G'$ satisfies $\theta_{i,n}\Leftrightarrow n_i(G')\geq n$.

\smallskip

One has $n_2(G)\geq n$ if and only if $G$ satisfies the following $\exists\forall$-sentence $\theta_{2,n}$: there exist $\ell$ finite subgroups $C_1,\ldots ,C_{\ell}$ of $G$ and, for every $1\leq i\leq \ell$, a finite subset $\lbrace g_{i,j}\rbrace_{1\leq j\leq n_i}$ of $N_G(C_i)$ such that:
\begin{itemize}
\item[$\bullet$]for every $g\in G$ and $1\leq i\neq j\leq n$, the groups $gC_ig^{-1}$ and $C_j$ are distinct;
\item[$\bullet$]the sum $n_1+\cdots+n_{\ell}$ is equal to $n$;
\item[$\bullet$]for every $1\leq i\leq \ell$, and for every $1\leq j\neq k\leq n_i$, the automorphisms $\mathrm{ad}(g_j)_{\vert C_i}$ and $\mathrm{ad}(g_k)_{\vert C_i}$ of $C_i$ are distinct.
\end{itemize}

\smallskip

One has $n_3(G)\geq n$ if and only if $G$ satisfies the following $\exists\forall$-sentence $\theta_{3,n}$: there exist $n$ finite subgroups $C_1,\ldots ,C_n$ of $G$ and $n$ elements $g_1\in N_G(C_1),\ldots,g_n\in N_G(C_n)$ of infinite order (i.e.\ satisfying $g_i^{K_G!}\neq 1$) such that:
\begin{itemize}
\item[$\bullet$]for every $g\in G$ and $1\leq i\neq j\leq n$, the groups $gC_ig^{-1}$ and $C_j$ are distinct;
\item[$\bullet$]for every $1\leq i\leq n$ and $g\in N_G(C_i)$, the subgroup $\langle g,g_i\rangle$ of $N_G(C_i)$ is virtually cyclic, i.e.\ $[g^{K_G!},g_i^{K_G!}]=1$.
\end{itemize}

\smallskip

One has $n_4(G)\geq n$ if and only if $G$ satisfies the following $\exists\forall$-sentence $\theta_{4,n}$: there exist $n$ finite subgroups $C_1,\ldots ,C_n$ of $G$ and, for every $1\leq i\leq n$, a couple of elements $(g_{i,1},g_{i,2})$ normalizing $C_i$ such that:
\begin{itemize}
\item[$\bullet$]for every $g\in G$ and $1\leq i\neq j\leq n$, the groups $gC_ig^{-1}$ and $C_j$ are distinct;
\item[$\bullet$]for every $1\leq i\leq n$, the subgroup $\langle g_{i,1},g_{i,2}\rangle$ is not virtually cyclic (i.e.\ $[g_{i,1}^{K_G!},g_{i,2}^{K_G!}]$ is non-trivial).
\end{itemize}

\smallskip

One has $n_5(G)\geq n$ if and only if $G$ satisfies the following $\exists\forall$-sentence $\theta_{5,n}$: there exist $2n$ finite subgroups $C_1,\ldots ,C_n$ and $C'_1\varsupsetneq C_1,\ldots,C'_n\varsupsetneq C_n$ of $G$ and, for every $1\leq i\leq n$, a couple of elements $(g_{i,1},g_{i,2})$ normalizing $C_i$, such that:
\begin{itemize}
\item[$\bullet$]for every $g\in G$ and $1\leq i\neq j\leq n$, the groups $gC_ig^{-1}$ and $C_j$ are distinct;
\item[$\bullet$]for every $1\leq i\leq n$, the subgroup $\langle g_{i,1},g_{i,2}\rangle$ is not virtually cyclic;
\item[$\bullet$]every element of $G$ that normalizes $C_i$ also normalizes $C'_i$.
\end{itemize}\end{proofL}

\section{Algorithm}

In this section, we shall prove the following theorem.

\begin{te}\label{algoalgoalgoalgo}There is an algorithm that, given a finite presentation of a virtually free group $G$ and a finite subset $X \subset G$, outputs `Yes' if the subgroup of $G$ generated by $X$ is $\exists\forall\exists$-elementary, and `No' otherwise.\end{te}

We shall use the following fact.

\begin{lemme}\label{lemmecomp}A subgroup $H$ of $G$ is $\exists\forall\exists$-elementary if and only if the three conditions of Proposition \ref{concon} are satisfied.\end{lemme}

\begin{proof}
If the conditions of Proposition \ref{concon} are satisfied, then either $H=G$, or $H$ is a proper subgroup and $G$ is a multiple legal large extension of $H$, by Proposition \ref{concon}. In both cases, the subgroup $H$ is $\exists\forall\exists$-elementary by Theorem \ref{legalteplus}. Conversely, if $H$ is $\exists\forall\exists$-elementary, then either $H=G$ or $H$ is a proper subgroup of $G$ and $G$ is a multiple legal large extension of $H$, by Theorem \ref{1}.\end{proof}

The proof of Theorem \ref{algoalgoalgoalgo} consists in showing that the conditions of Proposition \ref{concon} can be decided by an algorithm.

\subsection{Algorithmic tools} First, we collect several algorithms that will be useful in the proof of Theorem \ref{algoalgoalgoalgo}. 

\subsubsection{Solving equations in hyperbolic groups}The following theorem is the main result of \cite{DG10}.

\begin{te}\label{DG10}There exists an algorithm that takes as input a finite presentation of a hyperbolic group $G$ and a finite system of equations and inequations with constants in $G$, and decides whether there exists a solution or not.
\end{te}

\subsubsection{Computing a finite presentation of a subgroup given by generators}The following result is a particular case of Theorem 20 in \cite{BM16}.

\begin{te}\label{compu}There is an algorithm that, given a finite presentation of a hyperbolic and locally quasiconvex group $G$, and a finite subset $X$ of $G$, produces a finite presentation for the subgroup of $G$ generated by $X$.
\end{te}

Recall that a group is said to be \emph{locally quasiconvex} if every finitely generated subgroup is quasiconvex. Marshall Hall Jr.\ proved in \cite{Mar49} that every finitely generated subgroup of a finitely generated free group is a free factor in a finite-index subgroup, which shows in particular that finitely generated free groups are locally quasiconvex. It follows easily that finitely generated virtually free groups are locally quasiconvex. Thus, Theorem \ref{compu} applies when $G$ is virtually free.

\subsubsection{Basic algorithms}

\begin{lemme}\label{algofini}There is an algorithm that takes as input a finite presentation of a hyperbolic group and computes a list of representatives of the conjugacy classes of finite subgroups in this hyperbolic group.\end{lemme}

\begin{proof}There exists an algorithm that computes, given a finite presentation $\langle S \ \vert \ R\rangle$ of a hyperbolic group $G$, a hyperbolicity constant $\delta$ of $G$ (see \cite{Pa96}). In addition, it is well-known that the ball of radius $100\delta$ in $G$ contains at least one representative of each conjugacy class of finite subgroups of $G$ (see \cite{Bra00}). Moreover, two finite subgroups $C_1$ and $C_2$ of $G$ are conjugate if and only if there exists an element $g$ whose length is bounded by a constant depending only on $\delta$ and on the size of the generating set $S$ of $G$, such that $C_2=gC_1g^{-1}$ (see \cite{BH05}).\end{proof}

\begin{lemme}[\cite{DG11}, Lemma 2.5]\label{algon}There is an algorithm that computes a set of generators of the normalizer of any given finite subgroup in a hyperbolic group.\end{lemme}

\begin{lemme}[\cite{DG11}, Lemma 2.8]\label{algon2}There is an algorithm that decides, given a finite set $S$ in a hyperbolic group, whether $\langle S\rangle$ is finite, virtually cyclic infinite, or non virtually cyclic (finite or infinite).\end{lemme}

\begin{lemme}\label{algofini2}There is an algorithm that takes as input a finite presentation of a hyperbolic group $G$ and a finite subgroup $C$ of $G$ such that $N_G(C)$ is non virtually cyclic (finite or infinite), and decides whether or not $E_G(N_G(C))=C$.\end{lemme}

\begin{proof}By Lemma \ref{algofini}, one can compute some representatives $A_1,\ldots,A_k$ of the conjugacy classes of finite subgroups of $G$. Given an element $g\in G$, let $\theta_g(x)$ be a quantifier-free formula expressing the following fact: there exists an integer $1\leq i\leq k$ such that the finite set $\lbrace C,g\rbrace$ is contained in $xA_ix^{-1}$. Note that the group $\langle C,g\rangle$ is finite if and only if the existential sentence $\exists x \ \theta_g(x)$ is true in $G$.

One can compute a finite generating set $S$ for $N_G(C)$ using Lemma \ref{algon}. By Theorem \ref{DG10} above, one can decide if the following existential sentence with constants in $G$ is satisfied by $G$: there exist two elements $g$ and $g'$ such that
\begin{enumerate}
\item $g$ does not belong to $C$;
\item $\theta_g(g')$ is satisfied by $G$ (hence, the subgroup $C':=\langle C,g\rangle$ is finite);
\item for every $s\in S$, one has $sC's^{-1}=C'$.
\end{enumerate} 
Note that such an element $g$ exists if and only if $C$ is strictly contained in $E_G(N_G(C))$. This concludes the proof of the lemma.\end{proof}

The following lemma is an immediate corollary of Lemmas \ref{algofini}, \ref{algon}, \ref{algon2} and \ref{algofini2} above.

\begin{lemme}\label{algo5}There is an algorithm that takes as input a finite presentation of a hyperbolic group $G$ and computes the five numbers $n_1(G),\ldots,n_5(G)$ (see Definition \ref{5nombres}).
\end{lemme}

\subsection{Decidability of the first condition of Proposition \ref{concon}}

\begin{lemme}\label{condition1}There is an algorithm that, given a finite presentation of a virtually free group $G$ and a finite subset $X \subset G$ generating a subgroup $H=\langle X\rangle$, outputs `Yes' if $n_i(H)=n_i(G)$ for all $i\in \lbrace 1,2,3,4,5\rbrace$ and `No' otherwise.\end{lemme}

\begin{proof}
By Theorem \ref{compu}, there is an algorithm that takes as input a finite presentation $G=\langle S_G \ \vert \ R_G\rangle$ and $X$, and produces a finite presentation $\langle S_H \ \vert \ R_H\rangle$ for $H$. By Lemma \ref{algo5}, one can compute $n_i(G)$ and $n_i(H)$ for every $i\in \lbrace 1,2,3,4,5\rbrace$.\end{proof}

\subsection{Decidability of the second condition of Proposition \ref{concon}}

\begin{lemme}\label{condition2}There is an algorithm that, given a finite presentation of a virtually free group $G$ and a finite subset $X \subset G$ generating a subgroup $H=\langle X\rangle$, outputs `Yes' if $H$ is infinite and coincides with the one-ended factor of $G$ relative to $H$ (well-defined since $H$ is infinite), and `No' otherwise.\end{lemme}

\begin{proof}By Lemma \ref{algon2}, one can decide if $H$ is finite or infinite. By Lemma 8.7 in \cite{DG11}, one can compute a Stallings splitting of $G$ relative to $H$. Let $T$ be the Bass-Serre tree of this splitting. Let $U$ be the one-ended factor of $G$ relative to $H$ and let $u$ be the vertex of $T$ fixed by $U$. By Corollary 8.3 in \cite{DG11}, one can decide if there exists an automorphism $\varphi$ of $G$ such that $\varphi(H)=U$, which is equivalent to deciding if $U=H$. Indeed, if $\varphi(H)=U$, then $H$ fixes the vertex $u$ for the action of $G$ on $T$ twisted by $\varphi$. Thus, by definition of $U$ as the one-ended factor relative to $H$, the pair $(U,H)$ acts trivially on the tree $T$ for the action twisted by $\varphi$. Consequently, $\varphi(U)$ fixes $u$ as well. Therefore, one has $\varphi(U)=U=\varphi(H)$, and it follows that $U=H$ since $\varphi$ is an automorphism of $G$.\end{proof}

\subsection{Decidability of the third condition of Proposition \ref{concon}}

\begin{lemme}\label{condition3}There is an algorithm that, given a finite presentation of a virtually free group $G$ and a finite subset $X \subset G$ generating a subgroup $H=\langle X\rangle$, decides whether or not every finite subgroup of $G$ is conjugate to a subgroup of $H$.\end{lemme}

\begin{proof}
By Theorem \ref{compu}, there is an algorithm that takes as input a finite presentation $G=\langle S_G \ \vert \ R_G\rangle$ and $X$, and produces a finite presentation $\langle S_H \ \vert \ R_H\rangle$ for $H$. By Lemma \ref{algofini}, there is an algorithm that computes two lists $\lbrace A_1,\ldots ,A_n\rbrace$ and $\lbrace B_1,\ldots,B_n\rbrace$ of representatives of the conjugacy classes of finite subgroups of $G$ and $H$ respectively. Then, for every finite subgroup $A_i$ of $G$ in the first list, deciding if $A_i$ is conjugate in $G$ to $B_j$ for some $j\in \lbrace 1,\ldots ,n\rbrace$ is equivalent to solving the following finite disjunction of systems of equations with constants in $G$, which can be done using Theorem \ref{DG10}:
\[\theta(x) : \exists x \ (xA_ix^{-1}=B_1) \vee \ldots \vee (xA_ix^{-1}=B_n).\]
Hence, there is an algorithm that outputs `Yes' if every finite subgroup of $G$ is conjugate to a subgroup of $H$, and `No' otherwise.\end{proof}

\subsection{Proof of Theorem \ref{algoalgoalgoalgo}}

Theorem \ref{algoalgoalgoalgo} is an immediate consequence of Lemma \ref{lemmecomp} combined with Lemmas \ref{condition1}, \ref{condition2} and \ref{condition3}.

\section{$\exists^{+}$-elementary morphisms}

We prove Theorem \ref{44}.

\begin{te}Let $G$ be a group. Suppose that $G$ is finitely presented and Hopfian, or finitely generated and equationally noetherian. Then, every $\exists^{+}$-endomorphism of $G$ is an automorphism.\end{te}

\begin{proof}Let $\langle g_1,\ldots ,g_n\ \vert \ R(g_1,\ldots,g_n)=1\rangle$ be a presentation of $G$, with $R$ eventually infinite. Let $\varphi: G \rightarrow G$ be an $\exists^{+}$-endomorphism. For every integer $1\leq i\leq n$, let $h_i=\varphi(g_i)$. 

Note that there is a one-to-one correspondence between the set of homomorphisms $\mathrm{Hom}(G,G)$ and the set of solutions in $G^n$ of the system of equations $R(x_1,\ldots,x_n)=1$. If $G$ is finitely presentable, one can assume without loss of generality that the system of equations $R$ is finite. If $G$ is equationally noetherian, there is a finite subsystem $R_i(x_1,\ldots,x_n)=1$ of $R(x_1,\ldots,x_n)=1$ such that the sets $\mathrm{Hom}(G,G)$ and $\mathrm{Hom}(G_i,G)$ are in bijection, where $G_i$ denotes the finitely presented group $\langle g_1,\ldots ,g_n\ \vert \ R_i(g_1,\ldots,g_n)=1\rangle$. Hence, one can always assume without loss of generality that the system $R$ is finite.

Every element $h_i$ can be written as a word $w_i(g_1,\ldots,g_n)$, and the group $G$ satisfies the following existential positive formula:\[\mu(h_1,\ldots,h_n):\exists x_1 \ldots \exists x_n \ R(x_1,\ldots,x_n)=1 \wedge h_i=w_i(x_1,\ldots,x_n).\]Indeed, one can just take $x_i=g_i$ for every $1\leq i\leq n$, which shows that the formula $\mu(h_1,\ldots,h_n)$ is satisfied by $G$. Since this formula is existential positive and since the morphism $\varphi$ is $\exists^{+}$-elementary and $h_i=\varphi(g_i)$, the statement $\mu(g_1,\ldots,g_n)$ is true in $G$ too. As a consequence, there exist some elements $k_1,\ldots,k_n$ in $G$ such that $R(k_1,\ldots,k_n)=1$ and $g_i=w_i(k_1,\ldots,k_n)$ for every $1\leq i\leq n$. Let us define an endomorphism $\psi$ of $G$ by $\psi(g_i)=k_i$ for every $1\leq i\leq n$. 

Recall that $h_i=w_i(g_1,\ldots,g_n)$. Thus, one has \[\psi(h_i)=w_i(\psi(g_1),\ldots,\psi(g_n))=w_i(k_1,\ldots,k_n)=g_i.\]As a consequence, the composition $\psi\circ \varphi$ maps $g_i$ to itself, i.e.\ is the identity of $G$. It follows that $\psi$ is surjective. 

Last, recall that equationally noetherian groups are Hopfian. It follows that $\psi$ is an automorphism of $G$. Hence, $\varphi$ is an automorphism of $G$.\end{proof}

\renewcommand{\refname}{Bibliography}
\bibliographystyle{plain}
\bibliography{biblio}

\end{document}